\documentclass[12pt,a4paper,reqno]
{amsart}
\usepackage{amsfonts, amsmath, amssymb, amsthm, amsbsy, amscd, euscript}

\usepackage[english]{babel}

\usepackage[top=2.25cm, bottom=2.25cm, left=2.25cm, right=2.25cm]{geometry}

\usepackage{stmaryrd}

\usepackage{cmap}

\usepackage{tikz,tikz-cd}
\usetikzlibrary{arrows,positioning,automata,shadows,fit,shapes}

\usepackage[pdftex,
unicode=true,
colorlinks=false
]{hyperref}

\usepackage{DLdef1}

\renewcommand {\ssbegin}[2][*]
 {\refstepcounter{subsection}%
\if#1*
\addcontentsline{toc}{subsection}{\thesubsection.\hskip 1pc #2}%
\else
\addcontentsline{toc}{subsection}{\thesubsection.\hskip 1pc #2. #1}%
\fi
 \def \secno {\gdef \secno {}{\ssecfont
\thesubsection.\hskip 2ex}%
 }%
 \begin{#2}}

\renewcommand {\sssbegin}[2][*]
  {\refstepcounter{subsubsection}
\if#1*
\addcontentsline{toc}{subsubsection}{\thesubsubsection.\hskip 1pc #2}%
\else
\addcontentsline{toc}{subsubsection}{\thesubsubsection.\hskip 1pc #2. #1}
\fi
  \def \secno {\gdef \secno {}{\ssecfont \thesubsubsection.\hskip 2ex}%
  }%
   \begin{#2}}

\renewcommand {\parbegin}[2][*]
  {\refstepcounter{paragraph}
\if#1*
\addcontentsline{toc}{paragraph}{\theparagraph.\hskip 1pc #2}%
\else
\addcontentsline{toc}{paragraph}{\theparagraph.\hskip 1pc #2. #1}
\fi
  \def \secno {\gdef \secno {}{\ssecfont \theparagraph.\hskip 2ex}%
  }%
   \begin{#2}}

\renewcommand{\ssecfont}{\normalfont}

\newcommand{\dd}{\mathrm{d}}

\newcommand{\del}{\partial}
\newcommand{\what}{\widehat}

\DeclareMathOperator{\tsigma}{\tilde{\sigma}}

\newcommand{\trr}{\triangleright}

 \DeclareMathOperator{\im}{im}

\DeclareMathOperator{\Cl}{Cl}
\newcommand{\twedge}{{\textstyle{\bigwedge}}}

\title{Quantised $\fsl_2$-differential algebras}

\author{Andrey Krutov}
\address[Krutov]{Mathematical Institute of Charles University, Sokolovsk\'a 83, Prague, Czech Republic}
\email{andrey.krutov@matfyz.cuni.cz}

\author{Pavle Pand\v{z}i\'{c}}
\address[Pand\v zi\'c]{Department of Mathematics, Faculty of Science, University of Zagreb, Bijeni\v{c}ka cesta 30, 10000 Zagreb, Croatia}
\email{pandzic@math.hr}



\thanks{
  A.~Krutov was supported by the GA\v{C}R project 24-10887S
  and by HORIZON-MSCA-2022-SE-01-01 CaLIGOLA.
  P.~Pand\v{z}i\'{c} was supported by the project ``Implementation of cutting-edge research and its application as part of the Scientific Center of Excellence for Quantum and Complex Systems, and Representations of Lie Algebras'', PK.1.1.02, European Union, European Regional Development Fund.
  This article is based upon work from COST Action CaLISTA CA21109 supported by COST (European
  Cooperation in Science and Technology, \url{www.cost.eu}).
}

\begin{document}

\begin{abstract}
  We propose a~definition of a~quantised $\fsl_2$-differential algebra and show that
  the quantised exterior algebra (defined by Berenstein and Zwicknagl)  and the quantised Clifford algebra (defined by
  the authors) of~$\fsl_2$ are natural examples of such algebras.
\end{abstract}

\maketitle

\section{Introduction}
Let $\fg$ be a~Lie algebra.
H.~Cartan introduced the notion of $\fg$-differential algebras as a~generalisation
of algebras of differential forms on manifolds with $\fg$-action,  \cite{Cartan1951trans,HCartan1951}.
Later $\fg$-differential algebras appeared in the study of equivariant cohomology~\cite{GuilleminSternbergBook,AlekseevMeinrenken2000},
in Chern--\/Weil theory~\cite{AlekseevMeinrenken2005,MeinrenkenBook}, and in relation to (algebraic) Dirac
operators and Vogan's conjecture~\cite{AlekseevMeinrenken2000,HP1,HP2}.

We develop a new approach to quantisation of the notion of $\fg$-differential algebras which assumes that the quantum exterior algebra has classical
dimension. 
This is different from previous attempts to generalise the notion of $\fg$-differential algebras to the setting of
quantum groups and noncommutative geometry; for example,
see~\cite{AschieriCastellani1993,SchuppWattsZumino1993cartan,AschieriSchupp1996}.
Some of these works start with a~bicovariant calculus on a~quantum group, which usually does not have classical dimension. This
is in particular the case for $U_q(\fsl_2)$, see~\cite{WoroDc}, see also~\cite{JurcoCalculi} for the general case.
Other works assume the setting of triangular Hopf algebras which is not applicable to the setting of~$U_q(\fsl_2)$
since it is only quasitriangular, see~\cite{DrinfeldICM}.

In this paper
we propose a~definition of quantised $\fsl_2$-differential algebras based on the 3-dimensional quantised
adjoint $U_q(\fsl_2)$-module and give first examples,
certain quantised Clifford and exterior algebras. 
The advantage of our approach is that we start with the quantum exterior algebra defined by Berenstein and
Zwicknagl~\cite{BerensteinZwicknagl2008} of the classical dimension instead of a~bicovariant calculus.
We use the coboundary structure on the category of $U_q(\fsl_2)$-modules, see~\cite{DrinfeldQuasiHopf}. (As it
was shown in~\cite{HenriquesKamnitzer2006} such coboundary structure is related to the category of crystals.)

The paper is organised as follows. In \S\ref{sec:Pre} we recall necessary facts about the Drinfeld--\/Jimbo
quantum group~$U_q(\fsl_2)$, the quantised adjoint representation and its quantum exterior algebra.
In~\S\ref{sec:Clq} we recall the definition of the $q$-deformed Clifford algebra of~$\fsl_2$ introduced in~\cite{cubicDiracUqSL2} and
define Lie derivatives, contraction operators and the differential on it. We show that the defined operations
enjoy many features of their classical counterparts, in particular, Cartan's magic formula holds for them.
In~\S\ref{sec:Gen} we propose a~definition of a~quantised $\fsl_2$-differential algebra and show that the
quantised exterior and Clifford algebras of~$\fsl_2$ are examples of such algebras.

\subsection*{Acknowledgments}
We would like to thank P.~Aschieri, D.~Bashkirov, and R.~\'O~Buachalla for stimulating discussions.
The authors thank the organisers of the conference
``Quantum Groups and Noncommutative Geometry in Prague''
(11th to 15th September 2023, Prague, Czech Republic) for warm hospitality and financial support.
We are also grateful to the referee for useful suggestions.

\section{Preliminaries}\label{sec:Pre}
\subsection{$\fg$-differential algebras}
Let $\fg$ be a complex Lie algebra. Let  $\twedge[\xi]$ be the Grassmann algebra with
one generator~$\xi$, and let $\dd:=\del_\xi\in\Der\twedge[\xi]$ be the derivation with respect to~$\xi$.
Set $\what{\fg} :=  \fg\otimes\twedge[\xi] \inplus \Cee\dd$. Then
$\what{\fg} = \what{\fg}_{-1}\oplus\what{\fg}_{0}\oplus\what{\fg}_1$ is a $\Zee$-graded Lie superalgebra where
\[
  \what{\fg}_{-1} = \fg\otimes\xi,\qquad
  \what{\fg}_0 = \fg\otimes 1,\qquad
  \what{\fg}_1 = \Cee \dd.
\]
For $x\in\fg$, let  $L_x:= x \otimes 1\in\what{\fg}_{0}$, $\iota_x := x\otimes\xi \in \what{\fg}_{-1}$.
The non-zero bracket relations in~$\what{\fg}$ are defined as
\begin{gather}
  [L_x,\iota_y] = \iota_{[x,y]},\qquad
  [L_x,L_y] = L_{[x,y]},\qquad
  [\iota_x,\dd] = L_x
  \qquad\text{for all $x,y\in\fg$.}
\end{gather}

\subsubsection{Digression: semisimple Lie superalgebras}
Assume that $\fg$ is simple.
Let $\twedge(n)$ denote the Grassmann algebra with $n$ generators
$\xi_1,\ldots,\xi_n$. Then $\twedge(n)$ has a natural $\Zee$-grading given by~$\deg
\xi_i=1$.   Let $\fvect(0|n) := \Der\twedge(n)$. Clearly, $\fvect(0|n)$ is a $\Zee$-graded Lie superalgebra where $\deg \del_{\xi_i}=-1$. Let
$\fvect(0|n)_{-1}$ denotes the homogeneous component of degree~$-1$.
As it was shown in~\cite{Cheng1995}, 
 any semisimple Lie superalgebra is the direct sum of the following summands
\[
  \tilde{\mathfrak{s}}\otimes\twedge(n) \inplus \mathfrak{v},
\]
where $\mathfrak{s}$ is a simple Lie superalgebra, $\mathfrak{s}\subseteq\tilde{\mathfrak{s}}\subseteq\Der\mathfrak{s}$, and $\mathfrak{v}\subset\fvect(0|n)$ is such that
the projection
$\mathfrak{v}\to\fvect(0|n)_{-1}$ is onto.
In our case  (for $\what\fg$) we have that $n=1$, $\mathfrak{v}=\Span_\Cee(\del_{\xi})$, $\tilde{\mathfrak{s}}=\mathfrak{s}=\fg$.

\subsection{$\fg$-differential spaces and algebras}
A \emph{$\fg$-differential space} is a superspace~$V$, together with a $\what\fg$-module
structure~$\rho\colon\what\fg\to\End(V)$.
A \emph{$\fg$-differential algebra} is a superalgebra~$A$, equipped with a structure of $\fg$-differential space
such that $\rho(x)\in\Der A$ for all $x\in\what{\fg}$.
Observe that if $A$ is a $\fg$-differential algebra then the contraction operators~$\iota$ define a $\fg$-equivariant
representation of~$U(\what{\fg}_{-1}) \cong \twedge\fg$ on~$A$, where $U(\what{\fg}_{-1})$ is the universal
enveloping algebra of the Lie superalgebra~$\what{\fg}_{-1}$.
The idea of a~$\fg$-differential algebra is due to H.~Cartan~\cite{HCartan1951,Cartan1951trans}. We follow the terminology and
notation from~\cite{MeinrenkenBook}.

\sssbegin{Example}\label{ex:Ext}
Take $A = \twedge\fg^*$, equipped with the coadjoint action of~$\fg$
denoted by $L_x$ for $x\in\fg$. For $x\in\fg$ and $f\in\fg^\ast = \twedge^1\fg^\ast$
define the contraction operator by $\iota_xf = f(x)$. The odd map $\iota_x$ is extended
to~$\twedge\fg^\ast$ by the super Leibniz rule. Let $e_a$ be a~basis of~$\fg$ and $f_a$ be the corresponding
dual basis in~$\fg^\ast$. The Lie algebra differential on~$\twedge\fg^\ast$ may be written as
\[
  \dd_{\wedge}  = \frac12 \sum_a f_a\circ L_{e_a},
\]
with $f_a$ acting by the exterior multiplication.
Then $\twedge\fg^*$ is a $\fg$-differential algebra. One can show that
$H(\twedge\fg^*,\dd_\wedge)\cong (\twedge\fg^*)^\fg$.
\end{Example}

\sssbegin{Example}\label{ex:Cl}
Suppose that~$\fg$ has a nondegenerate invariant symmetric bilinear form~$B$ (for example, see review in~\cite{NIS}), used to identify $\fg\cong\fg^*$.  Let
$\Cl(\fg)$ be the Clifford algebra of~$\fg$ with respect to~$B$ defined by
\[
  \Cl(\fg) = T(\fg) / \left\langle x\otimes y + y\otimes x - 2B(x,y) \mid x,y\in\fg \right\rangle.
\]
Let $z_i$ be an~orthonormal basis of~$\fg$, then the Chevalley map (or quantisation)
$q_{\Cl}\colon \twedge(\fg)\to \Cl(\fg)$ is defined by
\[
  z_{i_1}\wedge \ldots \wedge z_{i_k} \mapsto z_{i_1}\ldots z_{i_k}\qquad\text{(and $1\mapsto1$)},
\]
where $1\leq i_1 < \ldots < i_{k}  \leq \dim\fg$.
Set
\[
  \gamma = - \frac{1}{12} \sum_{a,b,c=1}^{\dim\fg} B([z_a,z_b],z_c) z_a\wedge z_b\wedge z_c\in(\twedge^3\fg)^\fg.
\]
Define the map $\alpha\colon \fg\to\Cl(\fg)$ by
\[
  \alpha(x) = -\frac{1}{4}\sum_{a,b=1}^{\dim\fg} B(x,[z_a,z_b])z_az_b\quad\text{for $x\in\fg$}.
\]
The map~$\alpha$ extends to an~algebra homomorphism~$\alpha\colon U(\fg)\to\Cl(\fg)$.

The Clifford algebra~$\Cl(\fg)$ is a~filtered $\fg$-differential algebra with differential, Lie derivatives
and contractions given as
\[
  \dd_{\Cl} = [q_{\Cl}(\gamma), - ]_{\Cl},\qquad
  L_x = [\alpha(x), - ]_{\Cl},\qquad
  \iota_x = \frac12[x, - ]_{\Cl},\quad
  \text{for $x\in\fg$},
\]
where $[-,-]_{\Cl}$ denotes the supercommutator in~$\Cl(\fg)$.
The quantisation map~$q_{\Cl}\colon\twedge\fg\to\Cl(\fg)$ intertwines the Lie derivatives and contractions, but does
not intertwine the differential.
The cohomology of~$(\Cl(\fg),\dd_{\Cl})$ is trivial in all filtration degrees (except if~$\fg$ is abelian, in which case
$\dd_{\Cl}=0$); for example, see~\cite[\S7.1]{MeinrenkenBook}.
\end{Example}

\subsection{$U_q(\fsl_2)$}
Fix a~nonzero $q\in\Cee$ which is not a~root of unity.
The quantised enveloping algebra $U_q(\fsl_2)$ is the  associative algebra with unit
generated by the elements $E$, $F$, $K$, and~$K^{-1}$ subject to the relations
\begin{gather*}
  KE = q^2EK,\quad KF=q^{-2}FK,\quad KK^{-1}=K^{-1}K=1,\quad
  EF - FE = \frac{K-K^{-1}}{q-q^{-1}}.
\end{gather*}
A Hopf algebra structure on~$U_q(\fsl_2)$ is given by
\begin{gather*}
  \Delta E = E\otimes K + 1\otimes E,\quad
  \Delta F = F\otimes 1 + K^{-1}\otimes F,\quad
  \Delta K = K\otimes K,\quad
  \Delta K^{-1} = K^{-1}\otimes K^{-1},\\
  S(E) = - EK^{-1},\quad
  S(F) = - KF,\quad
  S(K^{-1}) = K,\quad
  S(K) = K^{-1},\\
  \eps(E) = \eps(F) = 0,\quad
  \eps(K) = \eps(K^{-1}) = 1,
\end{gather*}
where $\Delta$ is the coproduct, $S$ is the antipode, and $\eps$ is the counit.
In what follows we use Sweedler notation for the coproduct $\Delta x = \sum x_{(1)}\otimes x_{(2)}$.

Let $\fh$ be a~Cartan subalgebra of~$\fsl_2$, $\cP\subset\fh^\ast$ be the weight lattice of~$\fsl_2$,
and $\cP_{+}$ be the sublattice of dominant weights generated by the fundamental weight~$\pi$.
The category of finite dimensional type~1 modules over~$U_q(\fsl_2)$ is equivalent to the category of finite
dimensional $\fsl_2$ modules; for example, see~\cite[\S5.8]{EtingofTensorCat} or~\cite[\S3]{KSLeabh}. For $\lambda\in\cP_{+}$ we denote the corresponding
type~1 finite dimensional $U_q(\fsl_2)$-module with highest weight~$\lambda$ by~$V_\lambda$.

Let $\fsl_q(2)$ denote the vector subspace of~$U_q(\fsl_2)$ spanned by the elements
\[
  X = E,\qquad Z = q^{-2}EF - FE,\qquad Y =  KF.  
\]
The space~$\fsl_q(2)$ is closed with respect to the left adjoint action of $U_q(\fsl_2)$ on itself defined by
\[
  \ad_x y = \sum x_{(1)} y S(x_{(2)})\qquad\text{for $x,y\in U_q(\fsl_2)$.}
\]
It is easy to see that as a~$U_q(\fsl_2)$-module, $\fsl_q(2)$ is isomorphic to the quantised adjoint
representation~$V_{2\pi}$ of~$\fsl_2$.
In what follows we will use notation $\fsl_q(2)$ to emphasise that elements $X$, $Z$, and $Y$ belong to
$\fsl_q(2)\subset U_q(\fsl_2)$. In the case when $V_{2\pi}$ is treated as an abstract $U_q(\fsl_2)$-module and
in the case when we will construct quantum exterior and Clifford algebras, we will use the following notation
for basis elements in~$V_{2\pi}$:
\[
v_2 = X,\qquad v_0 = Z,\qquad  v_{-2} = Y.
\]

\subsection{Normalised braiding}
The following construction is due to Drinfeld~\cite{DrinfeldQuasiHopf}. Let $\mathsf{C}$ be a~braided monoidal
category linear over $\Cee[[\hbar]]$ and assume that the braiding satisfies
$\sigma_{W,V} \circ \sigma_{V,W} = \id_{V\otimes W} + O(\hbar)$. Then the map 
\[
  \tsigma_{V,W} = \sigma_{V,W} \circ (\sigma_{W,V}\circ \sigma_{V,W})^{-1/2},
\]
is called a normalised braiding and defines a~coboundary structure on~$\mathsf{C}$ in the sense of~\cite{DrinfeldQuasiHopf}.
For details
see~\cite[Exercise~8.3.25 on p.~202]{EtingofTensorCat}.
In particular, we have that $\tsigma^2=\id$.

The category of type one finite-dimensional $U_q(\fsl_2)$-module is a~braided monoidal category where the
braiding~$\sigma$ is given by the universal $R$-matrix; see~\cite[\S8.3]{EtingofTensorCat} for details.
The $R$-matrix braiding~$\sigma$ satisfies the above condition.
In what follows we denote by~$\tsigma$ the corresponding normalised braiding.

\subsection{Quantum exterior algebras}
Following~\cite{BerensteinZwicknagl2008} define the quantum exterior algebra~$\twedge_qV_{2\pi}$ of~$V_{2\pi}$
as 
\[
  \twedge_qV_{2\pi} = T(V_{2\pi}) /\langle v\otimes w + \tsigma(v\otimes w)\mid v,w\in V_{2\pi}\rangle
\]
The algebra~$\twedge_qV_{2\pi}$ is generated by $v_2$, $v_0$, $v_{-2}$ subject to the following relations
\begin{align*}
  v_2\wedge v_2 = {} &0,
  &v_{-2}\wedge v_{-2} = {}& 0,\\
  v_0\wedge v_2 = {}& {}-q^{-2} v_2\wedge v_0,
  &v_{-2}\wedge v_0 = {}& {}-q^{-2} v_0\wedge v_{-2},\\
  v_0\wedge v_0 = {}& \frac{(1-q^4)}{q^3} v_2\wedge v_{-2},
  &v_{-2}\wedge v_2 ={}&{}  -v_2\wedge v_{-2}.
\end{align*}
We note that $\twedge_qV_{2\pi}$ is a $\Zee$-graded super algebra in the braided monoidal category of type~1
finite-dimensional $U_q(\fsl_2)$-modules. The $\Zee_2$-grading corresponding to a~super algebra structure is
given by setting $p(v_2)=p(v_0)=p(v_{-2})=\od$, where $p(v)\in\Zee_2=\{\ev,\od\}$ denotes the parity of the element~$v$.
The algebra~$\twedge_qV_{2\pi}$ is (super)commutative with respect to normalised braiding, i.e.
$v \wedge w  = (-1)^{p(v)p(w)}\wedge \circ \tsigma(v\otimes w)$ for all parity homogeneous $v,w\in\twedge_qV_{2\pi}$.

\section{Motivating example: $\Cl_q(\fsl_2)$}\label{sec:Clq}
Fix a~non-zero constant $c\in\Cee[q,q^{-1}]$.
First recall from~\cite[\S2.7]{cubicDiracUqSL2} that $V_{2\pi}$ admits a~nondegenerate $U_q(\fsl_2)$-invariant bilinear form given by
\[
  \langle v_2, v_{-2} \rangle = c,\qquad
  \langle v_0, v_0 \rangle = q^{-3}(1+q^2)c,\qquad
  \langle v_{-2}, v_{2} \rangle = c q^{-2},
\]
Note that the form $\langle\cdot,\cdot\rangle$ is symmetric with respect to the normalised braiding~$\tsigma$,
i.e., $\langle \cdot,\cdot\rangle = \langle\cdot,\cdot\rangle \circ \tsigma$.
The $q$-deformed Clifford algebra of~$\fsl_2$ was defined in~\cite[\S3]{cubicDiracUqSL2} as filtered
deformation of~$\twedge_qV_{2\pi}$ by the bilinear form~$\langle\cdot,\cdot\rangle$:
\[
  \Cl_q(\fsl_2) = T(V_{2\pi}) /\langle v\otimes w + \tsigma(v\otimes w) - 2\langle v,w \rangle\mid v,w \in V_{2\pi}\rangle,
\]
As it was shown in~\cite[Lemma~3.3]{cubicDiracUqSL2}, the algebra $\Cl_q(\fsl_2)$ is generated by~$v_2$,
$v_0$, $v_{-2}$ satisfying the following relations
\begin{align*}
  v_2 v_2 = {} &0,
  &v_{-2} v_{-2} = {}& 0,\\
  v_0 v_2 = {}& {}-q^{-2} v_2 v_0,
  &v_{-2} v_0 = {}& {}-q^{-2} v_0 v_{-2},\\
  v_0 v_0 = {}& \frac{1-q^4}{q^3} v_2 v_{-2} + \frac{q^2+1}{q} c 1,
  &v_{-2} v_2 ={}&{}  -v_2 v_{-2} + \frac{q^2+1}{q^2} c 1.
\end{align*}
It is easy to see that $\Cl_q(\fsl_2)$ is a~filtered super algebra in the (braided) monoidal category of
$U_q(\fsl_2)$-modules. We note that $\Cl_q(\fsl_2)$ is a~filtered $U_q(\fsl_2)$-module where the
elements of~$U_q(\fsl_2)$ act by operators of degree~$0$.

The quantum exterior algebra~$\twedge_qV_{2\pi}$ and the $q$-deformed Clifford algebra~$\Cl_q(\fsl_2)$ can be
considered as results of the deformation quantisation of the classical exterior algebra~$\twedge\fsl_2$. The
corresponding Poisson structures on~$\twedge\fsl_2$ are studied in~\cite{MomentMapsRmat}.

\subsection{$\tsigma$-commutators}
For $x,y\in\Cl_q(\fsl_2)$ homogeneous with respect to parity  set
\[
  [x,y]_{\tsigma} := \left(m_{\Cl_q} - (-1)^{p(x)p(y)}m_{\Cl_q}\circ\tsigma\right)(x\otimes y),
\]
where $m_{\Cl_q}$ denotes the multiplication map in~$\Cl_q(\fsl_2)$.
The map $[-,-]_{\tsigma}$ is extended to~$\Cl_q(\fsl_2)$ by linearity.
By construction $[-,-]_{\tsigma}$ is $U_q(\fsl_2)$-equivariant since it is composed from equivariant maps.

\sssbegin{Lemma}\label{lem:tCommSkew}
The bracket $[-,-]_{\tsigma}$ is $\tsigma$-skew-symmetric:
\begin{gather*}
  [\omega,\mu]_{\tsigma} = - (-1)^{p(\omega)p(\mu)} [-,-]_{\tsigma} \circ \tsigma (\omega\otimes\mu)\qquad
  \text{for $\omega,\mu\in\Cl_q(\fsl_2)$,}
\end{gather*}
and has the filtration degree~$-1$.
\end{Lemma}
\begin{proof} The $\tsigma$-skew-symmetricity follows form the definition of~$[-,-]_{\tsigma}$.  
  Since $\Cl_q(\fsl_2)$ is a~filtered deformation of the $\tsigma$-supercommutative
  algebra~$\twedge_qV_{2\pi}$, it follows from the definition that the bracket $[-,-]_{\tsigma}$ has the
  filtration degree~$-1$.
\end{proof}

\subsection{The $\alpha$ and $\beta$ maps and Lie derivatives}
The quantum moment map (in the sense of~\cite{Lu1993}) is the algebra map $\alpha_q\colon U_q(\fsl_2)\to
\Cl_q(\fsl_2)$ defined in~\cite[\S3.7]{cubicDiracUqSL2}. It is given on generators by
\begin{gather*}
  \alpha_q(E) = - \frac{q}{(1+q^2)c}v_2v_0,\qquad
  \alpha_q(F) = - \frac{q^2}{(1+q^2)c}v_0v_{-2},\\
  \alpha_q(K) = \frac{q^3-q}{(1+q^2)c}v_2v_{-2} + q^{-1},\qquad
  \alpha_q(K^{-1}) = -\frac{q^3-q}{(1+q^2)c}v_2v_{-2} + q.
\end{gather*}
Since $\alpha_q$ is an~algebra map it follows that
\[
\alpha_q(Z) = \frac1cv_2v_{-2} - 1,\qquad
  \alpha_q(Y) = -\frac{q}{(1+q^2)c}v_0v_{-2}.
\]
As it was shown in~\cite[Lemma~3.7.1]{cubicDiracUqSL2} the inner $U_q(\fsl_2)$-action defined by~$\alpha_q$
coincides with the natural one:
\[
  x \trr \omega = \sum \alpha_q(x_{(1)}) \omega \alpha_q(S(x_{(2)}))
  \qquad\text{for $x\in U_q(\fsl_2)$, $\omega\in\Cl_q(\fsl_2)$},
\]
where $x \trr \omega$ denotes the $U_q(\fsl_2)$-action on~$\Cl_q(\fsl_2)$.

Following the classical situation we define Lie derivatives on~$\Cl_q(\fsl_2)$ with respect to elements
of~$U_q(\fsl_2)$ by
\[
  L_x \omega := x \trr \omega\qquad
  \text{for $x\in U_q(\fsl_2)$, $\omega\in\Cl_q(\fsl_2)$}.
\]

Classically, the map~$\alpha$ defined the (adjoint) action of~$\fg$ by taking the (super)commutator; see Example~\ref{ex:Cl}.
This is no longer true in the quantum case for $\tsigma$-commutators.
Define a~linear map $\beta_q \colon \fsl_q(2) \to \Cl_q(\fsl_2)$ by
\[
  \beta_q(X) =  \frac{1+q^2}{q}\alpha_q(X),\qquad
  \beta_q(Y) =  \frac{1+q^2}{q}\alpha_q(Y),\qquad
  \beta_q(Z) =  \frac{1+q^2}{q}\alpha_q(Z).
\]
The definition of~$\beta_q$ is motivated by the following lemma.

\sssbegin{Proposition}
The $\beta_q$-map defines the quantum Hamiltonian with respect to the $\tsigma$-commu\-tator for the action of
elements of $\fsl_q(2)\subset U_q(\fsl_2)$ on~$\Cl_q(\fsl_2)$. Namely, we have that
\[
  L_x \omega = [\beta_q(x), \omega]_{\tsigma}\qquad
  \text{for $x\in\fsl_q(2)$, $\omega\in\Cl_q(\fsl_2)$.}
\]
\end{Proposition}
\begin{proof}
  For $X\in\fsl_q(2)$ and $v_2\in\Cl_q(\fsl_2)$, we have that
  \begin{align*}
    [\beta_q(X), v_2]_{\tsigma} = {}
    & {} - \frac{1}{2c}v_2v_0 v_2 + \frac{1}{2c}v_2v_2v_0
      = \frac{1}{2q^{2}c}v_2v_2v_0 = 0 = X\trr v_2.
  \end{align*}
  The computations for other elements of~$\fsl_q(2)$ and~$\Cl_q(\fsl_2)$ are analogous.
\end{proof}

\sssbegin{Remark}
One can check that $\beta_q$ defines a~morphism of ``quantum brackets'' in the following sense:
$[\beta_q(x),\beta_q(y)]_{\tsigma} = \beta_q(\ad_xy)$ for all $x,y\in\fsl_q(2)$.
Moreover, it follows from Lemma~\ref{lem:tCommSkew} that $\ad \circ \tsigma = - \ad$ on $\fsl_q(2)$.
\end{Remark}

\subsection{Differential}
Recall from~\cite[\S3.4]{cubicDiracUqSL2} that the element
\[
  \gamma_q = - \frac{1}{2c^2}(c v_0 + v_2v_0v_{-2}) \in \Cl_q^{(3)}(\fsl_2)
\]
squares to a~scalar. Therefore, we can define a~differential on~$\Cl_q(\fsl_2)$ by
\[
  \dd_{\Cl_q}\omega_q = [\gamma_q, \omega]_{\tsigma}
  = \gamma_q \omega - (-1)^{p(\omega)} \omega \gamma_q,
  \qquad \omega\in\Cl_q(\fsl_2).
\]

\sssbegin{Proposition}
We have that
\begin{enumerate}
\item $\dd_{\Cl_q} x = 2\beta_q(x)$ for $x \in \fsl_q(2)$;

\item the differential $\dd_{\Cl_q}$ is $U_q(\fsl_2)$-equivariant.
\end{enumerate}
\end{Proposition}
\begin{proof}
  (1) We have that
  \begin{align*}
    \dd_{\Cl_q} v_2 = {}
    &{} - \frac{1}{2c^2}( v_2(c v_0 + v_2v_0v_{-2}) + (c v_0 + v_2v_0v_{-2})v_2)\\
    {}={}& - \frac{1}{2c^2}( c v_2v_0 - q^{-2} c v_2v_0 - v_2v_0v_2v_{-2} + \frac{q^2+1}{q^2}c v_2v_0)\\
    {}={}&  - \frac{1}{2c^2} \frac{cq^2 - c + q^2c + c}{q^2} v_2v_2
           = - \frac{1}{c} v_2v_0 = 2\beta_q(X),\\
    \dd_{\Cl_q} v_0 = {}
    &{} - \frac{1}{2c^2}( v_0(c v_0 + v_2v_0v_{-2}) + (c v_0 + v_2v_0v_{-2})v_0)\\
    {}={}& -\frac{1}{2c^2}\left( \frac{2c(1-q^4)}{q^3} v_2v_{-2} + \frac{2c^2(q^2+1)}{q} - 2q^{-2}v_2v_0v_0v_{-2}\right)\\
    {}={}& -\frac{1}{2c^2}\left( \frac{2c(1-q^4)}{q^3} v_2v_{-2} + \frac{2c^2(q^2+1)}{q} - \frac{2c(q^2+1)}{q^{3}}v_2v_{-2}\right)\\
    {}={}& \frac{q^2+q}{2cq}(v_2v_{-2} - c) = 2\beta_q(Z),\\
    \dd_{\Cl_q} v_0 = {}
    &{} - \frac{1}{2c^2}( v_{-2}(c v_0 + v_2v_0v_{-2}) + (c v_0 + v_2v_0v_{-2})v_{-2})\\
    {}={}& - \frac{1}{2c^2}(-q^{-2}c v_0v_{-2} + q^{-2} v_{2}v_{-2}v_{-2}v_0
           + \frac{q^2+1}{q^2}cv_0v_{-2} + cv_0v_{-2})\\
    {}={}& - \frac{1}{c}v_0v_{-2} = 2\beta_q(Y).
  \end{align*}
  
  (2) The equivariance of~$\dd_{\Cl_q}$ follows from the equivariance of the bracket~$[-,-]_{\tsigma}$.
\end{proof}

\sssbegin{Remark}
We note that the differential is dual to the quantised adjoint action of $\fsl_q(2)$ on itself (a~version of
``quantum Lie bracket'') in the following sense: $ \iota_x\iota_y \dd_{\Cl_q} z = \langle \ad_xy,z \rangle$ for $x,y,z\in\fsl_q(2)$.
\end{Remark}

\subsection{Contractions}
Following the classical case, see Example~\ref{ex:Cl}, define 
\[
  \iota_x  \omega = \frac12[ x, \omega]_{\tsigma},\qquad x\in V_{2\pi}, \omega \in \Cl_q(\fsl_2).
\]
This definition is motivated by the fact that for linear $v\in V_{2\pi} \subset \Cl^{(1)}_q(\fsl_2)$ we have that
\[
  \iota_x v = \frac12[x,v]_{\tsigma}  =\frac12( xv + m_{\Cl_q}\circ\tsigma(x\otimes v))
  = \frac12( 2 \langle x,v \rangle) = \langle x, v \rangle,
\]
where $m_{\Cl_q}$ is the multiplication map in~$\Cl_q(\fsl_2)$.
Furthermore, $\iota_x$ has the filtration degree~$-1$.

\sssbegin{Proposition}
For $x,y\in V_{2\pi}$ let $x_i,y_i\in V_{2\pi}$ be defined by
$\tsigma(x\otimes y) = \sum_i y_i\otimes x_i$, then for all $\omega\in\Cl_q(\fsl_2)$
the contraction operators satisfy
\[
  \iota_x\iota_y\omega + \sum_i\iota_{y_i}\iota_{x_i}\omega = 0.
\]
Hence the map $\iota\colon V_{2\pi} \to \End(\Cl_q(\fsl_2))$ extends to a~$U_q(\fsl_2$)-equivariant morphism of superalgebras
$\iota\colon \twedge_qV_{2\pi} \to \End(\Cl_q(\fsl_2))$.
\end{Proposition}
\begin{proof}
  For $v_{2},v_{0}\in V_{2\pi}=\twedge^1_qV_{2\pi}$ we have that
  \[
    v_2\otimes v_0 + \tsigma(v_2\otimes v_0) = \frac{2}{1+q^4}\left(q^2 v_0\otimes v_2+ v_2\otimes v_0\right).
  \]
  For $v_2v_0v_{-2}\in\Cl_q(\fsl_2)$ we have that
  \begin{align*}
    \iota_{v_2}\iota_{v_0}(v_2v_0v_{-2}) =
    {} & \iota_{v_2}\left((q^2-1)(q^2+1)q^{-3}c^2 - (1+q^2)q^{-1}v_2v_{-2}\right)
         = (1+q^2)q^{-1}c^2 v_2,\\
    \intertext{and}
    \iota_{v_0}\iota_{v_2}(v_2v_0v_{-2}) =
    {} & \iota_{v_0}(cv_2v_0) = - (1+q^2)q^{-3}c^2 v_2 = - q^{-2} \iota_{v_2}\iota_{v_0}(v_2v_0v_{-2}).
  \end{align*}
  The computations for other elements of~$\twedge_qV_{2\pi}$ and~$\Cl_q(\fsl_2)$ are analogous.
  The $U_q(\fsl_2)$-equivariance follows from the equivariance of the bracket $[-,-]_{\tsigma}$.
\end{proof}

\sssbegin{Remark}\label{rem:DerIota}
It is subject of further research to describe the derivation property for the action of $\twedge_q V_{2\pi}$
on~$\Cl_q(\fsl_2)$. We show below that several obvious attempts to ensure this property do not work.

We consider the following possible approach: let us assume that for $x\in\fsl_q(2)$ and $\omega,\mu\in\Cl_q(\fsl_2)$, 
\[
  \iota_x (\omega \mu) = \iota_x(\omega) \mu + (-1)^{p(\omega)}\sum_i\omega_i\iota_{x_i}(\mu),
\]
where $\omega_i$ and $x_i$ are defined by the following options
\begin{enumerate}
\item $\sigma(x\otimes \omega) = \sum \omega_i \otimes x_i$, \label{it:Rmat}
\item $\sigma^{-1}(x\otimes \omega) = \sum \omega_i \otimes x_i$,\label{it:RmatInv}
\item $\tsigma(x\otimes \omega) = \sum \omega_i \otimes x_i$,\label{it:RmatNorm}
\end{enumerate}
In what follows we show that none of these options works in the example of $x=v_0$, $\omega = v_2$,
and $\mu =  v_{-2}$. First note that
\[
  \iota_{v_0} (v_2v_{-2}) = \frac{c(1-q^2)}{q^2}v_0,\qquad
  \iota_{v_0}(v_2) v_{-2} = 0.
\]

\underline{Case~\ref{it:Rmat}.} We have that
\[
  \sigma(v_0\otimes v_2) = v_2\otimes v_0 + q^{-2}(q^4-1)v_0\otimes v_2.
\]
Therefore,
\[
  \iota_x(\omega) \mu + (-1)^{p(\omega)}\sum_i\omega_i\iota_{x_i}(\mu) =
  - v_2\iota_{v_0}(v_{-2}) - \frac{q^4-1}{q^2}v_0 \iota_{v_2}v_{-2}
  = - \frac{c(q^4-1)}{q^2}v_0 \neq \iota_{v_0}(v_2v_{-2}).
\]

\underline{Case~\ref{it:RmatInv}.} We have that
\[
  \sigma^{-1}(v_0\otimes v_2) = v_2\otimes v_0.
\]
Therefore,
\[
  \iota_x(\omega) \mu + (-1)^{p(\omega)}\sum_i\omega_i\iota_{x_i}(\mu) =
  -  v_2\iota_{v_0}(v_{-2}) = 0 \neq \iota_{v_0}(v_2v_{-2}).
\]

\underline{Case~\ref{it:RmatNorm}.} We have that
\[
  \tsigma(v_0\otimes v_2) = \frac{2q^2}{1+q^4}v_2\otimes v_0 + \frac{q^4-1}{1+q^4}v_0\otimes v_2.
\]
Therefore,
\begin{align*}
  \iota_x(\omega) \mu + (-1)^{p(\omega)}\sum_i\omega_i\iota_{x_i}(\mu) 
  ={}&{}- \frac{2q^2}{1+q^4}v_2\iota_{v_0}(v_{-2}) - \frac{q^4-1}{1+q^4}v_0 \iota_{v_2}v_{-2}\\
  {}={}&{} - \frac{c(q^4-1)}{1+q^4}v_0 \neq \iota_{v_0}(v_2v_{-2}).
\end{align*}

Perhaps one could use approach similar to~\cite[\S6]{AschieriCastellani1993} or consider a~braided Hopf
algebra structure on~$\twedge_qV_{2\pi}$ similar, for example, to one constructed~\cite{NicholsGrass}.
\end{Remark}

In view of the difficulties explained in Remark \ref{rem:DerIota}, it is surprising that we have the following theorem.

\ssbegin{Theorem}
For $x\in\fsl_q(2)= V_{2\pi}$ the operators $L_x$, $\iota_x$, and $\dd_{\Cl}$ on~$\Cl_q(\fsl_2)$ satisfy Cartan's magic formula
\[
  L_x = \iota_x \circ \dd_{\Cl_q} + \dd_{\Cl_q} \circ \iota_x.
\]
In particular, cochain maps~$L_x$ are homotopic to~0, with $\iota_x$ as homotopy operators. Therefore, $L_x$
induces the zero action on cohomology.
\end{Theorem}
\begin{proof}
  Direct computations.
  For example, for $v_2\in\fsl_q(2)= V_{2\pi}$
  and $v_{-2}\in\Cl_q(\fsl_2)$ we have that
  \[
    L_{v_2}v_{-2} = v_0
  \]
  and
  \[
    \iota_{v_2}\dd_{\Cl_q}v_{-2} + \dd_{\Cl_q}  \iota_{v_2} v_{-2} = {}
    -\frac{1}{c}\iota_{v_2}v_0v_{-2} + c\dd_{\Cl_q}(1)  = v_0.
  \]
  The computations for other cases are similar.
\end{proof}

\ssbegin{Remark}\label{rem:Clrho}
First note that the element~$\gamma_q$ generates a $U_q(\fsl_2)$-invariant subalgebra
in~$\Cl_q(\fsl_2)$. Moreover, the element~$\gamma_q$ satisfies
\( \gamma_q^2 = \frac{1+q^2}{4cq}  \).
Therefore, by the universal property of Clifford algebras, we have that
\[
  \Cl_q(\fsl_2)^{U_q(\fsl_2)} = \Cl(P_q(\fsl_2), B_q),
\]
where $P_q(\fsl_2)$ is the space of primitive invariants spanned by~$\gamma_q$ equipped with nondegenerate
symmetric bilinear form~$B_q$ given by $B_q(\gamma_q,\gamma_q) = \frac{1+q^2}{4cq}$.
It is now easy to see that $H(\Cl_q(\fsl_2),\dd_{\Cl_q}) = 0$.
  
Since the element~$\gamma_q$ is $U_q(\fsl_2)$-invariant we have that
$[\omega,\gamma_q]_{\tsigma} = -(-1)^{p(\omega)} [\gamma_q,\omega]_{\tsigma}$ for all parity homogeneous $\omega\in\Cl_q(\fsl_2)$.
Therefore, for $x\in\fsl_q(2)$ we have that
\[
  \iota_x\gamma = [\tfrac12 x,\gamma_q]_{\tsigma} 
  = [\gamma_q,\tfrac12 x]_{\tsigma} = \tfrac12 \dd_{\Cl_q}(x)  = \beta_q(x) \in \im\alpha_q.
\]
Moreover, direct computations show that for $x\in\fsl_q(2)$ we have
\[
  \iota_x\gamma_q\cdot \gamma_q^\ast  = x,
\]
where $\gamma_q^\ast = \frac{4q c}{1+q^2}\gamma_q$ is the dual to~$\gamma_q$ with respect to~$B_q$.
This leads to the quantum analogue of the $\rho$-decomposition from~\cite{KostantRho}:
\[
  \Cl_q(\fsl_2) = \Cl(P_q(\fsl_2), B_q) \otimes \im \alpha_q.
\]
We emphasise that in this case the braided tensor product of algebras in the braided monoidal category of
type~1 finite-dimensional $U_q(\fsl_2)$-modules reduces to the usual tensor product of algebras since the elements of
$\Cl(P_q(\fsl_2),B_q)$ are $U_q(\fsl_2)$-invariant.
\end{Remark}

\section{The general definition}\label{sec:Gen}
\ssbegin{Definition}
A supervector space~$W$ is called a~quantised $\fsl_2$-differential space if it is equipped with
\begin{enumerate}
\item Lie derivatives $L_x\in\End(W)$ for $x\in U_q(\fsl_2)$ which define a $U_q(\fsl_2)$-module
  structure on~$W$;
\item a~$U_q(\fsl_2)$-equivariant action $\iota\colon \twedge_qV_{2\pi}\otimes W\to W$ of $\twedge_qV_{2\pi}$;
\item a~$U_q(\fsl_2)$-equivariant differential $\dd_W\colon W \to W$;
\item such that they satisfy  Cartan's magic formula
\[
  L_x = \iota_x \circ \dd_{W} + \dd_{W} \circ \iota_x\qquad\text{for $x\in\fsl_q(2)$}.
\]

\end{enumerate}
A morphism between two quantised $\fsl_2$-differential spaces is a~morphism in the category of
$U_q(\fsl_2)$-modules which intertwines contractions and differentials (and also Lie derivatives).
\end{Definition}

\ssbegin{Definition}
An algebra~$A$ is called a quantised $\fsl_2$-differential algebra
if it is a~quantised $\fsl_2$-differential space such that
\begin{enumerate}
\item the Lie derivatives satisfy
  \[
    L_x (ab ) =\sum (L_{x_{(1)}}a)(L_{x_{(2)}}b)\qquad
    \text{for $a,b\in A$, $x\in U_q(\fsl_2)$},
  \]
  in other words, $A$ is an~algebra in the monoidal category of $U_q(\fsl_2)$-modules;
\item the differential~$\dd_A$ satisfies the (graded) Leibniz rule.
\end{enumerate}
A morphism between two quantised $\fsl_2$-differential algebras is an~algebra morphism in the category of
$U_q(\fsl_2)$-modules which intertwines contractions and differentials (and also Lie derivatives).
\end{Definition}

It is subject of further research to describe the derivation property for the action of $\twedge_q V_{2\pi}$
on a~quantised $\fsl_2$-differential algebra; see the discussion in Remark~\ref{rem:DerIota}.

\subsection{Quantum exterior algebra}
In this subsection we show that the quantum exterior algebra of the 3-dimensional quantised adjoint
$U_q(\fsl_2)$-module~$V_{2\pi}$ is a~quantised $\fsl_2$-differential algebra in the sense of our definition.

First note that the associated graded algebra of~$\Cl_q(\fsl_2)$ is the quantum exterior
algebra~$\twedge_qV_{2\pi}$. For $x\in\fsl_q(2)$, the associated graded maps
$L_x\colon\twedge_qV_{2\pi}\to V_{2\pi}$ to the Lie derivatives
$L_x\colon\Cl_q(\fsl_2)\to \Cl_q(\fsl_2)$ define an action of $U_q(\fsl_2)$ since the filtration
on~$\Cl_q(\fsl_2)$ is compatible with $U_q(\fsl_2)$-action.

The differential $\dd_{\Cl_q}$ has filtered degree one. Therefore, we
can define the associated graded map~$\dd_{\wedge_q}\colon \twedge_qV_{2\pi}\to \twedge_qV_{2\pi}$ which is
$U_q(\fsl_2)$-equivariant by construction.
It is easy to see, c.f. Example~\ref{ex:Ext}, that $\dd_{\wedge_q}$ is nonzero only for
\[
  \dd_{\wedge_q}(v_2) = - \frac{1}{c}v_2\wedge v_0,\quad
  \dd_{\wedge_q}(v_0) = \frac{1+q^2}{q c} v_2\wedge v_{-2},\quad
  \dd_{\wedge_q}(v_{-2}) = - \frac{1}{c} v_0\wedge v_{-2}.
\]
It is straightforward to check that $\dd_{\wedge_q}^2=0$ and that it satisfies the graded Leibniz rule, so it
defines a differential on~$\twedge_qV_{2\pi}$. Moreover,  we have the quantised version of the
formula for differential, see Example~\ref{ex:Ext},
\[
  \dd_{\wedge_q} = \frac{q^2}{1+q^4}\left(
    \frac{1}{c}v_{-2}L_X + \frac{q^3}{(1+q^2)c}v_0 L_Z + \frac{q^2}{q^2}v_{2}L_Y
    \right).
\]
Note that the formulas for the differential depend on the parameter $c$ since we identify
$\twedge_qV_{2\pi}^\ast$ with~$\twedge_qV_{2\pi}$ via the bilinear form $\langle\cdot,\cdot\rangle$.

Similarly, for $x\in\twedge_qV_{2\pi}$ we define the contraction operator
$\iota_x\colon \twedge_qV_{2\pi}\to\twedge V_{2\pi}$ as the associated graded map for the contraction
operator on~$\Cl_q(\fsl_2)$. By construction, the operators $\iota_x$ define a $U_q(\fsl_2)$-equivariant representation of
$\twedge_qV_{2\pi}$. In particular, we have that
\begin{align*}
  & \iota_{v_2}v_2 = 0,
  & & \iota_{v_0}v_2 = 0,
  & &\iota_{v_{-2}}v_2 = q^{-2}c, \\
  & \iota_{v_2}v_0 = 0,
  & &\iota_{v_0}v_0 = q^{-3}(1+q^2)c,
  & &\iota_{v_{-2}}v_0 = 0, \\
  &\iota_{v_2}v_{-2} = c,  
  & & \iota_{v_0}v_2 = 0,
  & & \iota_{v_{-2}}v_{-2} = 0, \\
  & \iota_{v_2}v_2\wedge v_{0} = 0,
  & &\iota_{v_0}v_2\wedge v_0 =  -\tfrac{1+q^2}{q^3}cv_2,
  & &\iota_{v_{-2}}v_2\wedge v_0 = cv_0, \\  
  & \iota_{v_2}v_2\wedge v_{-2} = -cv_2,
  & &\iota_{v_0}v_2\wedge v_{-2} =  \tfrac{1-q^2}{q^2}cv_0,
  & & \iota_{v_{-2}}v_2\wedge v_{-2} = q^{-2}cv_{-2}, \\  
  &\iota_{v_2}v_{0}\wedge v_{-2} = -cv_0,
  & &\iota_{v_0}v_0\wedge v_{-2} =  \tfrac{1+q^2}{q}cv_{-2},
  & &\iota_{v_{-2}}v_0\wedge v_{-2} = 0, \\  
  &\iota_{v_2}v_2\wedge v_{0}\wedge v_{-2} = cv_2\wedge v_0,
  & &\iota_{v_0}v_2\wedge v_0\wedge v_{-2} =  - \tfrac{1+q^2}{q}cv_2\wedge v_{-2},
  & &\iota_{v_{-2}}v_2\wedge v_0\wedge v_{-2} =  cv_0\wedge v_{-2}.
\end{align*}

Cartan's magic formula for associated graded maps $L_x$, $\iota_x$, and $\dd_{\wedge_q}$
on~$\twedge_qV_{2\pi}$ follows from the fact that Cartan's magic formula on~$\Cl_q(\fsl_2)$ has degree~$0$.
It can also be checked by direct computations as follows.
First note that for elements of $\twedge^0V_{2\pi}$ and $\twedge^3V_{2\pi}$ Cartan's magic formula holds trivially.
The operator $\dd_{\wedge_q}\circ \iota_x$ acts by zero on~$\twedge^1V_{2\pi}$.
We have that
\[
  \iota_{v_{2}}\dd_{\wedge_q} v_{-2} =-\frac{1}{c} \iota_{v_2}( v_0\wedge v_{-2}) = v_0 = L_{v_2}v_{-2}.
\]
The operator $\iota_x\circ\dd_{\wedge_q}$ acts by zero on~$\twedge^2V_{2\pi}$.
We have that
\[
  \dd_{\wedge_q}\iota_{v_2}(v_0\wedge v_{-2}) =-  c \dd_{\wedge_q}v_0 = -\frac{1+q^2}{q} v_2\wedge v_{-2}
  = L_{v_2}(v_0\wedge v_{-2}).
\]
Therefore, we have proved the following theorem.

\ssbegin{Theorem}
The algebras $\Cl_q(\fsl_2)$ and $\twedge_qV_{2\pi}$ are quantised $\fsl_2$-differential algebras.
\end{Theorem}

\ssbegin{Remark}
Similarly to the case of~$\Cl_q(\fsl_2)$, see Remark~\ref{rem:Clrho}, we can now compute the cohomology 
of~$\twedge_qV_{2\pi}$ using Cartan's magic formula.
First note that the subalgebra of $U_q(\fsl_2)$-invariant elements in~$\twedge_qV_{2\pi}$ is spanned by
$1$, $v_2\wedge v_0\wedge v_{-2}$. Therefore, we have the quantised analogue of Hopf--\/Koszul--\/Samelson theorem
\[
  H(\twedge_qV_{2\pi},\dd_{\wedge_q}) = (\twedge_qV_{2\pi})^{U_q(\fsl_2)}
  = \twedge P_{\wedge_q}(\fsl_2),
\]
where $P_{\wedge_q}(\fsl_2)$ is the space of primitive invariants spanned by~$v_2\wedge v_0 \wedge v_{-2}$.

We note that 
\[
  \iota_{v_2\wedge v_0 \wedge v_{-2}}(v_2\wedge v_0 \wedge v_{-2}) = \iota_{v_2\wedge v_0}(c v_0 \wedge v_{-2})
  = c^2\tfrac{1+q^2}{q^2} \iota_{v_2}v_{-2} =  c^3\tfrac{1+q^2}{q^2}.
\]
Therefore, we can recover the form~$B_q$ from Remark~\ref{rem:Clrho} via contraction operator just as in the
classical case of~\cite{KostantRho}.
\end{Remark}

\providecommand{\bysame}{\leavevmode\hbox to3em{\hrulefill}\thinspace}
\providecommand{\MR}{\relax\ifhmode\unskip\space\fi MR }
\providecommand{\MRhref}[2]{%
  \href{http://www.ams.org/mathscinet-getitem?mr=#1}{#2}
}
\providecommand{\href}[2]{#2}


\end{document}